\documentclass[11pt,letterpaper]{article}

\usepackage{mypreprint2010}
\usepackage{graphicx}
\usepackage{amssymb}

\newenvironment{myproof}{\begin{proof}}{\end{proof}}
\newenvironment{thm}{\begin{theorem}}{\end{theorem}}
\newenvironment{lem}{\begin{lemma}}{\end{lemma}}
\newenvironment{rmk}{\begin{remark}}{\end{remark}}

\usepackage{amsmath}

\usepackage{algorithm}
\usepackage{algpseudocode}

\def\bfA{\boldsymbol{A}} 
\def\bfB{\boldsymbol{B}} \def\bfb{\boldsymbol{b}}
\def\bfC{\boldsymbol{C}} 
\def\bfD{\boldsymbol{D}} 
\def\bfE{\boldsymbol{E}}

\def\bfI{\boldsymbol{I}}

\def\bfL{\boldsymbol{L}} 
\def\bfM{\boldsymbol{M}}

\def\bfP{\boldsymbol{P}} 
 
\def\bfR{\boldsymbol{R}} 
\def\bfS{\boldsymbol{S}} 
\def\bfT{\boldsymbol{T}} 
\def\bfU{\boldsymbol{U}} 
\def\bfV{\boldsymbol{V}} 
\def\bfW{\boldsymbol{W}} 
\def\bfX{\boldsymbol{X}} 
 \def\bfy{\boldsymbol{y}}
\def\bfZ{\boldsymbol{Z}} \def\bfz{\boldsymbol{z}}
\def\bfnull{\boldsymbol{0}}

\def\Bs{\boldsymbol{B}_{\bot}}

\def\Ph{\boldsymbol{\widehat{P}}}
\def\H2{\mathcal{H}_2}

\DeclareMathOperator*{\colsp}{span}
\DeclareMathOperator*{\real}{Re}
\DeclareMathOperator*{\imag}{Im}

\begin{document}

\date{\today}

\title{ADI iteration for Lyapunov equations: a tangential approach and adaptive shift selection}

\author{Thomas Wolf\footnotemark[2]\ , \ Heiko K. F. Panzer\footnotemark[2]\ \footnotemark[3]\ \ and Boris Lohmann\footnotemark[2]}

\maketitle

\renewcommand{\thefootnote}{\fnsymbol{footnote}}
\footnotetext[2]{Institute of Automatic Control,
				Technische Universit\"at M\"unchen,
				Boltzmannstr. 15, D-85748 Garching (\textbraceleft \texttt{thomas.wolf,panzer,lohmann}\textbraceright\texttt{@tum.de})}
\footnotetext[3]{Partially supported by the Cusanuswerk Scholarship Award Programme, which is gratefully acknowledged.}
\renewcommand{\thefootnote}{\arabic{footnote}}

\begin{abstract}
A new version of the alternating directions implicit (ADI) iteration for the solution of large-scale Lyapunov equations is introduced.
It generalizes the hitherto existing iteration, by incorporating tangential directions in the way they are already available for rational Krylov subspaces.
Additionally, first strategies to adaptively select shifts and tangential directions in each iteration are presented.
Numerical examples emphasize the potential of the new results.

\end{abstract}

\section{Introduction}
We consider the generalized algebraic Lyapunov equation
\begin{equation}\label{eq:lyap}
	\bfA\bfP\bfE^T + \bfE\bfP\bfA^T +\bfB\bfB^T = \bfnull,
\end{equation}
with $\bfA,\bfE,\bfP \in \mathbb{R}^{n \times n}$ and $\bfB \in \mathbb{R}^{n \times m}$.
The order $n$ is assumed to be large and $m$ small, $m \ll n$.
Furthermore, the matrix $\bfE$ has full rank, det$(\bfE) \ne 0$.
The standard formulation of \eqref{eq:lyap} results from setting $\bfE$ to identity.
The solution $\bfP$ of the Lyapunov equation \eqref{eq:lyap} represents the so-called controllability Gramian of a linear time-invariant dynamical system.
It is used e.\,g. in model reduction by \emph{truncated balanced realization} (TBR), see \cite{antoulas} for details.

Direct (also called dense) numerical solvers of \eqref{eq:lyap} are available in matrix computation software \cite{direct_lyap}.
However, for high-dimensional problems with large $n$, they become inappropriate due to high memory requirements and execution time.
A remedy is to employ iterative methods that can exploit the sparsity of the matrices $\bfA$, $\bfE$ and $\bfB$.
These methods typically compute a low-rank factor $\bfZ$, having less columns than rows.
The approximation $\Ph \approx \bfP$ is then given by the low-rank formulation $\Ph := \bfZ\bfZ^H$.
This reduces memory requirements, because only the low-rank factor $\bfZ$ has to be stored, instead of the full $n$-by-$n$ matrix $\Ph$.

Here we focus on the \emph{alternating directions implicit} (ADI) iteration in its low-rank formulation \cite{Li_ADI,penzl_smith}.
It successively accumulates the low-rank factor $\bfZ$ by a block of $m$ columns.
This, however, can be a drawback if $m>1$ and a high number of iterations is necessary for sufficient approximation:
a low-rank factor $\bfZ$, growing by $m$ columns per iteration, might become too large for reasonable processing.

In this work we propose a \emph{tangential} generalization of the ADI iteration, for which the low-rank factor $\bfZ$ grows by only one column per iteration---irrespective of $m$.
The hitherto existing ADI iteration based on blocks is included as a special case.
The tangential ADI iteration is most appropriately performed in combination with an adaptive shift selection.
We also address this by suggesting a first procedure that determines new shifts and tangential directions on the fly.

In Section~\ref{sec:prelim} relevant preliminaries are reviewed.
A reformulation of the ADI iteration is presented in Section~\ref{sec:re_ADI}, from which the tangential version is derived in Section~\ref{sec:tangetial_ADI}.
The adaptive shift selection is addressed in Section~\ref{sec:adapt} and numerical examples are given in Section~\ref{sec:num_ex}.

Throughout the paper, we use the following notation:
$\bfI_k$ denotes the $k\times k$ identity matrix.
Complex conjugation of a matrix $\bfX$ is denoted by $\overline{\bfX}$; $\bfX^T$ means transposition and $\bfX^H = \overline{\bfX}^T$ is transposition with complex conjugation.
The real and imaginary part of a complex $\bfX$ are given by $\bfX = \real(\bfX) + \imath\imag(\bfX)$, where $\imath$ is the imaginary unit.
The open right half of the complex plane is denoted by $\mathbb{C}_+$ and
$\Lambda(\bfA)$ denotes the spectrum of a quadratic matrix $\bfA$.

\section{Preliminaries}\label{sec:prelim}
We review the low-rank ADI iteration for solving \eqref{eq:lyap} in the following.

\subsection{The ADI Iteration}
For a given set $\mathcal{S} = \{s_1,s_2, \ldots, s_k \}$ of complex valued shifts $s_i \in \mathbb{C}_+$, the blocks $\bfZ_i \in \mathbb{C}^{n \times m}$ of the low-rank factor $\bfZ = [\bfZ_1, \,\bfZ_2, \, \ldots, \, \bfZ_k]$ are given by the following ADI iteration for $i = 2,\ldots,k$, \cite{Li_ADI}:
\begin{equation}\label{eq:LRADI}
\begin{aligned} 
	\bfZ_1 \!&= \!\sqrt{2\, \real (s_1) }  \left( \bfA - s_1\bfE \right)^{-1} \bfB, \\
	\bfZ_{i} \!&= \!\sqrt{\frac{ \real (s_{i})}{ \real (s_{i-1}) } } \!
	\left[ \bfI \!+\! ( s_{i} \!+\! \overline{s}_{i-1} ) \! \left( \bfA \!-\! s_{i}\bfE \right)^{-1}\!\bfE \right]\! \bfZ_{i-1}. 
\end{aligned}
\end{equation}
As $\bfA,\bfE$ and $\bfB$ are real,
we assume the set $\mathcal{S}$ to be closed under conjugation, causing also $\Ph=\bfZ\bfZ^H$ to be real.
Typically, the set $\mathcal{S}$ is chosen a priori and reused periodically until convergence occurs \cite{Li_ADI,diss_saak}.
The residual $\bfR \in \mathbb{R}^{n \times n}$ is the remainder in the original Lyapunov equation \eqref{eq:lyap}, after substituting $\bfP$ by $\Ph$:
\begin{equation}\label{eq:res}
	\bfR := \bfA \Ph \bfE^T + \bfE \Ph \bfA^T + \bfB\bfB^T.
\end{equation}
In \cite{paed_MCMDS} and \cite{Wolf_ADI} it was independently shown, that the residual can be factorized as $\bfR = \Bs\Bs^T$---with real $\Bs \in \mathbb{R}^{n \times m}$, if the set $\mathcal{S}$ is closed under conjugation.
One way to compute $\Bs$ is
\begin{align}
	\Bs & = \bfB + \bfE\bfZ\bfL^T, \quad \textrm{with}  \label{eq:Bs} \\
	\bfL &:= [\, \sqrt{2\real(s_1)}\bfI_m, \, \ldots,\, \sqrt{2\real(s_k)}\bfI_m \; ].
\end{align}
The notation $\Bs$ stems from the fact, that a specific projection can be associated to the ADI iteration:
$\Bs$ is the residual after projecting $\bfB$, and it fulfills the Petrov-Galerkin condition, see \cite{Wolf_ADI} for details.
The low-rank formulation $\bfR = \Bs\Bs^T$ allows to compute the induced matrix $2$-norm of the residual $\|\bfR\|_2$ by the maximum eigenvalue of the $m$-by-$m$ matrix $\Bs^T\Bs$: $\|\bfR\|_2 \!=\! \max\, \Lambda(\Bs^T\Bs)$.
This can be done with low numerical effort, which is why one often uses $\|\bfR\|_2$ as a convergence criterion in the ADI iteration.

\subsection{The Problem}\label{sec:problem}
The low-rank factor $\bfZ$ of the ADI iteration \eqref{eq:LRADI} gains a new block of $m$ columns in each step.
For problems that require many iterations for convergence, i.\,e. until $\|\bfR\|_2$ is small enough, the final $\bfZ$ might become too large for reasonable processing.
It would be desirable, that for every shift $s_i$ only one column, or as many columns as absolutely necessary, are added to the low-rank factor $\bfZ$.
An ad hoc solution would be to substitute $\bfB$ in \eqref{eq:LRADI} by $\bfB\bfb$, where $\bfb \in \mathbb{C}^{m}$ is referred to as a \emph{tangential direction} in the following.
However, as every block/column $\bfZ_i$ in the iteration \eqref{eq:LRADI} originates from its predecessor $\bfZ_{i-1}$, this tangential direction is fixed for all time, and the left over directions in $\bfB$ would be completely neglected.
Obviously, this generally cannot result in a good approximation.

Our first contribution will be a re-formulation of the ADI iteration \eqref{eq:LRADI}.
This re-formulation is more natural, as it incorporates the low-rank factor $\Bs$ of the residual, and it is essential for our second contribution:
we show, how each shift $s_i$ can be equipped with an individual tangential direction $\bfb_i$.
We will denote this variation with \emph{tangential (low-rank) ADI iteration} (T-LR-ADI) in the following.
As it is advisable to employ T-LR-ADI with adaptively selected shifts $s_i$ and tangential directions $\bfb_i$, a first procedure of this kind is our third contribution.
Because the results require the recently established connection of the ADI iteration to Krylov subspace projection methods \cite{Wolf_ADI}, we briefly review this in the following subsections.

\subsection{Rational Krylov Subspace Method}\label{sec:RKSM}
Another iterative approach to solve \eqref{eq:lyap} is the \emph{rational Krylov subspace method} (RKSM) \cite{Druskin_2011}, which is based on projections by Krylov subspaces.
For the ease of presentation, we assume distinct shifts $s_i \ne s_j$, for $i \ne j$, in the set $\mathcal{S}$.
The associated rational block Krylov subspace $\mathcal{K}_{b}$ is defined as
\begin{equation}\label{eq:Krylov}
	\mathcal{K}_{b} := \colsp \left\{ (\bfA-s_1\bfE)^{-1}\bfB, \, \ldots, \, (\bfA-s_k\bfE)^{-1}\bfB \right\}.
\end{equation}
As the set $\mathcal{S}$ is assumed to be closed under conjugation, a real basis $\bfV \in \mathbb{R}^{n \times km}$ of the subspace \eqref{eq:Krylov} can be computed.
Typically, a Gram-Schmidt orthogonalization is employed, resulting in an orthogonal basis $\bfV^T\bfV = \bfI_{km}$.
Together with a matrix $\bfW \in \mathbb{R}^{n \times km}$, an oblique projection can be defined, and the respective projections of $\bfA,\bfE$ and $\bfB$ read as
\begin{equation}\label{eq:ArErBr}
	\bfA_k := \bfW^T\bfA\bfV, \quad
	\bfE_k := \bfW^T\bfE\bfV, \quad
	\bfB_k := \bfW^T\bfB.
\end{equation}
Solving the projected Lyapunov equation
\begin{equation}\label{eq:proj_lyap}
	\bfA_k\bfP_k\bfE_k^T + \bfE_k\bfP_k\bfA_k^T +\bfB_k\bfB_k^T = \bfnull,
\end{equation}
of dimension $km \times km$ by direct methods for $\bfP_k$ results in the approximation $\Ph_{\textrm{RKSM}} := \bfV\bfP_k\bfV^T$, which fulfills the \emph{Petrov-Galerkin} condition on the residual: $\bfW^T\bfR\bfW=\bfnull$.
Typically, one chooses $\bfW := \bfV$, see \cite{Druskin_2011} for details.
A low-rank formulation of the residual was introduced in \cite{Wolf_2013_ACC}.

\subsection{Relation of ADI Iteration to Krylov Subspace Methods}\label{sec:Krylov_ADI}
Li and White \cite{Li_ADI} showed that the ADI basis $\bfZ$ spans the rational block Krylov subspace: $\colsp(\bfZ) = \colsp(\bfV) = \mathcal{K}_b$.
First links between the ADI iteration and RKSM were introduced in \cite{Druskin_Convergence} and \cite{flagg_ADI} for particular sets $\mathcal{S}$.
Yet the full connection was recently established in \cite{Wolf_ADI}.
As our contributions are based on this connection, we briefly review it in the following.

The remaining degree of freedom in RKSM is the choice of $\bfW$, which defines the direction of projection.
It was shown that there exists a $\bfW$, such that the approximations of RKSM and ADI coincide: $\Ph_{\textrm{RKSM}} \!=\! \Ph$.
In other words, the approximation $\Ph_{\textrm{RKSM}}$ of RKSM is equal to the approximation $\Ph$ of the ADI iteration \eqref{eq:LRADI}---for a particular choice of $\bfW$.
This means that there are two possible ways to compute the ADI approximation for an arbitrary given set of shifts $\mathcal{S}$: either employ the iteration \eqref{eq:LRADI}, or set up the rational Krylov subspace, perform the particular projection and solve the reduced Lyapunov equation \eqref{eq:proj_lyap}.

In order to carry out the second approach to compute the ADI approximation,
one possible way is to calculate the desired direction of projection, i.\,e. $\bfW$.
However, it was shown in \cite{Wolf_ADI}, that the explicit formulation of $\bfW$ can be avoided.
This relies on the following fact:
the approximations of RKSM and ADI coincide, $\Ph_{\textrm{RKSM}} \!=\! \Ph$, if and only if the projected matrices \eqref{eq:ArErBr} fulfill $\Lambda(-\bfE_k^{-1}\bfA_k) \!=\! \mathcal{S}$.
That means that the Ritz values, mirrored along the imaginary axis, correspond to the set of shifts in the rational Krylov subspace.
An algorithm that exploits this fact, is given in the following theorem.
\begin{thm}[\cite{Wolf_ADI}]\label{thm:ADI_Krylov}
	For an arbitrary given set of shifts $\mathcal{S} = \{s_1,s_2, \ldots, s_k \}$, the approximation $\Ph=\bfZ\bfZ^H$ of the ADI iteration \eqref{eq:LRADI} can be equivalently computed by the following three steps:
	\begin{description}
		\item[Step 1:] Set up the basis $ \bfV= \left[ (\bfA-s_1\bfE)^{-1}\bfB, \, \ldots, \, (\bfA-s_k\bfE)^{-1}\bfB \right]$ of the Krylov subspace \eqref{eq:Krylov}, together with the matrices
		\begin{equation}\label{eq:example1}
			\bfS = \left[ \begin{array}{ccc} s_1\bfI_m \\ & \ddots \\ & & s_k\bfI_m \end{array} \right]
			\; \textrm{and} \;\;
			\bfL = \left[ \bfI_m, \, \ldots, \, \bfI_m \right].
		\end{equation}
		\item[Step 2:] Solve the reduced Lyapunov equation
		\begin{equation}\label{eq:X}
			\bfS\bfX + \bfX\bfS^H - \bfL^H\bfL = \bfnull,
		\end{equation}
		 for $\bfX$ by direct methods.
		\item[Step 3:] Compute the approximation $\Ph = \bfV\bfX^{-1}\bfV^H$.
	\end{description}
\end{thm}
\begin{rmk}
	The approach of Theorem~\ref{thm:ADI_Krylov} therefore also avoids the explicit construction of the projected matrices \eqref{eq:ArErBr}.
	The solution $\bfP_k$ of the reduced Lyapunov solution \eqref{eq:proj_lyap} is given by the inverse $\bfX^{-1} =: \bfP_k$.
	For details on the link between RKSM and ADI please refer to \cite{Wolf_ADI}.
\end{rmk}
\begin{rmk}
	For the ease of presentation, we chose a particular basis $\bfV$ of the Krylov subspace in Theorem~\ref{thm:ADI_Krylov}.
	It would be also possible to perform the three steps for an arbitrary basis $\widetilde{\bfV}$ with $\colsp(\widetilde{\bfV})=\colsp(\bfV)$.
	Then there exists a matrix $\bfT \in  \mathbb{C}^{km \times km}$, such that $\widetilde{\bfV} \!=\! \bfV\bfT$.
	The matrices $\bfS$ and $\bfL$ would become $\widetilde{\bfS} \!=\! \bfT^{-1}\bfS\bfT$ and $\widetilde{\bfL} \!=\! \bfL\bfT$.
	This is due to the fact, that every basis of a Krylov subspace can equivalently be formulated as the solution of a Sylvester equation \cite{Gallivan_2004,Wolf_2012_MathMod},
	\begin{equation}\label{eq:syl}
		\bfA\bfV - \bfE\bfV\bfS = \bfB\bfL,
	\end{equation}
	where the eigenvalues of $\bfS \in \mathbb{C}^{km \times km}$ correspond to the set of shifts $\mathcal{S} = \Lambda(\bfS)$.
	This means that for a given basis $\bfV$, the corresponding $\bfS$ and $\bfL$ solve \eqref{eq:syl}.
	Please note, that all matrices $\bfV,\bfS$ and $\bfL$ can be computed by an Arnoldi-like approach.
	Therefore, a different basis $\bfV$ of the Krylov subspace would only affect the first step in Theorem~\ref{thm:ADI_Krylov}, the steps 2 and 3 can be performed without modification.
\end{rmk}
\begin{rmk}
	If $\Ph$ of the ADI iteration is computed by Theorem~\ref{thm:ADI_Krylov}, the residual $\bfR$ can also be given in a low-rank formulation.
	However, it is slightly different from the formulation in the ADI iteration \eqref{eq:Bs} and given in \cite{Wolf_ADI} by
	\begin{equation}\label{eq:Bs_Krylov}
		\Bs = \bfB + \bfE\bfV\bfX^{-1}\bfL^H.
	\end{equation}
	As $\bfL$ can be complex for arbitrary bases $\bfV$, $\bfL^T$ of \eqref{eq:Bs} changes to $\bfL^H$ here.
\end{rmk}

\section{New Formulation of the ADI iteration}\label{sec:re_ADI}
The Krylov-based approach by Theorem~\ref{thm:ADI_Krylov} can be readily extended to incorporate tangential directions:
for a given set of shifts $\{s_1,s_2, \ldots, s_k \}$, $s_i \in \mathbb{C}$ and tangential directions $\{\bfb_1,\bfb_2, \ldots, \bfb_k \}$, $\bfb_i \in \mathbb{C}^m$, the tangential rational Krylov subspace $\mathcal{K}_{t}$ is defined as
\begin{equation}\label{eq:Krylov_t}
	\mathcal{K}_{t} \!:=\! \colsp \left\{ (\bfA \!\!-\!\! s_1\bfE)^{-1}\bfB\bfb_1, \ldots, (\bfA \!\!-\!\! s_k\bfE)^{-1}\bfB\bfb_k \right\}.
\end{equation}
Substituting the block subspace $\mathcal{K}_{b}$ \eqref{eq:Krylov} in Theorem~\ref{thm:ADI_Krylov} with the tangential subspace $\mathcal{K}_{t}$ \eqref{eq:Krylov_t} already yields the desired result.
For example, choose the basis $\bfV \!=\! [ (\bfA-s_1\bfE)^{-1}\bfB\bfb_1,\ldots,$$(\bfA-s_k\bfE)^{-1}\bfB\bfb_k]$, then the Syl\-ves\-ter equation \eqref{eq:syl} holds for
\begin{equation}\label{eq:example2}
	\bfS = \textrm{diag} \left(  s_1,\, \ldots,\, s_k \right)
	\; \textrm{and} \;\;
	\bfL = \left[ \bfb_1, \, \ldots, \, \bfb_k \right].
\end{equation}
Taking this $\bfS$ and $\bfL$ to solve \eqref{eq:X} for $\bfX$, the desired approximation $\Ph\!=\!\bfV\bfX^{-1}\bfV^H$---which will be the goal of our tangential ADI iteration---is already found.
Its rank is then $k$ instead of $mk$.

This approach is a new result and it already solves the problem we formulated in Section~\ref{sec:problem}:
every shift we use in the ADI iteration can now be equipped with an individual tangential direction.
Although this specifies a reasonable method for solving Lyapunov equations, it has the drawback that in an iterative procedure, a reduced Lyapunov equation of increasing dimension has to be solved by direct methods (see step 2 of Theorem~\ref{thm:ADI_Krylov}).
This results in additional numerical effort, which is nonexistent in the original ADI iteration \eqref{eq:LRADI}.
Our main contribution, presented in Section~\ref{sec:tangetial_ADI}, is to set up an iterative algorithm in an ADI fashion similar to \eqref{eq:LRADI}, which computes the desired result---without requiring direct solves of reduced Lyapunov equations.
This algorithm is denoted as \emph{tangential low-rank alternating directions implicit} (T-LR-ADI) iteration, as already mentioned.
T-LR-ADI is based on the re-formulation of \eqref{eq:LRADI} by incorporating the residual factor $\Bs$, presented in the following lemma.
\begin{lem}\label{thm:BSADI}
	Define $\bfB_{\bot,0} := \bfB$, then the ADI iteration \eqref{eq:LRADI} is equivalent to the following iteration for $i = 1,\ldots,k$:
	\begin{equation}\label{eq:BSADI}
	\begin{aligned}
		\bfZ_i &= \sqrt{2\, \real (s_i) }  \left( \bfA - s_i\bfE \right)^{-1} \bfB_{\bot,i-1}, \\
		\bfB_{\bot,i} &= \bfB_{\bot,i-1} + \sqrt{2\, \real (s_i) }\bfE\bfZ_i.
	\end{aligned}
	\end{equation}
\end{lem}
\begin{myproof}
	Consider the first step of the ADI iteration \eqref{eq:LRADI}, with the approximation $\Ph_1 := \bfZ_1\bfZ_1^H$ and the residual\footnote{For complex $s_1$, the set $\mathcal{S}\!=\!\{s_1\}$ is not closed under conjugation. Then the residual factor $\Bs$ becomes complex and the residual reads as $\bfR\!=\!\Bs\Bs^H$ instead of $\bfR\!=\!\Bs\Bs^T$.} $\bfR_1 := \bfB_{\bot,1}\bfB_{\bot,1}^H$, such that
	\begin{equation}\label{eq:res2}
		\bfA \Ph_1 \bfE^T + \bfE \Ph_1 \bfA^T + \bfB\bfB^T = \bfB_{\bot,1}\bfB_{\bot,1}^H,
	\end{equation}
	with $\bfB_{\bot,1} = \bfB + \sqrt{2\, \real (s_1) } \bfE\bfZ_1$.
	Subtracting \eqref{eq:res2} from \eqref{eq:lyap} yields
	\begin{equation}\label{eq:lyap_bs}
		\bfA (\bfP - \Ph_1) \bfE^T + \bfE (\bfP - \Ph_1) \bfA^T + \bfB_{\bot,1}\bfB_{\bot,1}^H = \bfnull,
	\end{equation}
	which can be interpreted as follows:
	the error $\bfP\!-\!\Ph_1$ after the first step solves the new Lyapunov equation \eqref{eq:lyap_bs} with the right-hand side $\bfB_{\bot,1}\bfB_{\bot,1}^H$.
	Applying the first step of the ADI iteration \eqref{eq:LRADI} to \eqref{eq:lyap_bs}, results in
	\begin{equation}
		\bfZ_2 = \sqrt{2\, \real (s_2) }  \left( \bfA - s_2\bfE \right)^{-1} \bfB_{\bot,1},
	\end{equation}
	leading to the second approximation $\Ph_2 \!=\! \bfZ_2\bfZ_2^H$ and the second residual $\bfR_2 \!=\! \bfB_{\bot,2}\bfB_{\bot,2}^H$ with $\bfB_{\bot,2} \!=\! \bfB_{\bot,1} + \sqrt{2\, \real (s_2) } \bfE\bfZ_2$.
	The overall error $\bfP-\Ph_1-\Ph_2 \!=\! \bfP - [\bfZ_1,\bfZ_2][\bfZ_1,\bfZ_2]^H$ then solves the Lyapunov equation for the right-hand side $\bfB_{\bot,2}\bfB_{\bot,2}^H$.
	By repeating this procedure for $k$ steps, the final residual reads as
	\begin{equation}\label{eq:Bs_k}
		\bfB_{\bot,k} = \bfB + \sum_{i=1}^{k} \sqrt{2\, \real (s_i) } \bfE\bfZ_i.
	\end{equation}
	As $\bfB_{\bot,k}$ in \eqref{eq:Bs_k} is equal to $\Bs$ in \eqref{eq:Bs}, the approximations $\Ph$ of the ADI iteration \eqref{eq:LRADI} and the new formulation \eqref{eq:BSADI} have to be equal, which completes the proof.
\end{myproof}
\begin{rmk}
	The re-formulation \eqref{eq:BSADI} of the ADI iteration was independently found by Benner and K{\"u}rschner, see \cite{Paed_real_Syl} and \cite{Paed_ADI_Shifts}.
	The authors use it to suggest a new method to adaptively select shifts in the ADI iteration.
	We will compare T-LR-ADI to this approach by a numerical example in Section~\ref{sec:num_ex}.
\end{rmk}
Lemma~\ref{thm:BSADI} provides a more natural formulation of the ADI iteration than \eqref{eq:LRADI}, because it incorporates the residual factor $\Bs$ as an integral part of the iteration.
It therefore also updates the current residual on the way as a byproduct.
Surprisingly, this works with unchanged numerical effort:
the main effort in both formulations \eqref{eq:LRADI} and \eqref{eq:BSADI} is to solve a linear set of equations for every iterate $\bfZ_i$.
The remaining operations in both iterations are a single matrix-vector product with $\bfE$ and a weighted sum of two $n$--by--$m$ blocks.

\section{Tangential Low-Rank ADI}\label{sec:tangetial_ADI}
In the previous formulation \eqref{eq:LRADI}, every iterate $\bfZ_i$ originates from its predecessor $\bfZ_{i-1}$.
The main advantage of the new formulation \eqref{eq:BSADI} is that each $\bfZ_i$ originates from the recent residual factor $\bfB_{\bot,i-1}$ instead.
Following the proof of Lemma~\ref{thm:BSADI}, each $\bfB_{\bot,i}$ can be interpreted as the new right-hand side in the Lyapunov equation after $i$ steps of the ADI iteration.
In other words: the algorithm is restarted with $\bfB_{\bot,i}$.
This allows to use an individual tangential direction $\bfb_i \in \mathbb{C}^m$ for every iterate by substituting $\bfB_{\bot,i}$ in \eqref{eq:BSADI} with $\bfB_{\bot,i}\bfb_i$.
It remains to show, that this new iteration results in the same approximation $\Ph$ as the Krylov-based approach introduced in Theorem~\ref{thm:ADI_Krylov} with $\mathcal{K}_t$ instead of $\mathcal{K}_b$ (as mentioned in the beginning of Section~\ref{sec:re_ADI}), and how the residual factor $\bfB_{\bot,i}$ changes by incorporating tangential directions.
Different from the block iteration \eqref{eq:BSADI}, real and complex conjugated shifts have to be distinguished in T-LR-ADI.
Both cases are treated separately in the two following subsections.

\subsection{Tangential ADI for Real Shifts}
\begin{thm}\label{thm:Tang_ADI_real}
	Define $\bfB_{\bot,0} := \bfB$, and assume exclusively real shifts $s_i \!\in\! \mathbb{R}$ and real tangential directions $\bfb_i \!\in\! \mathbb{R}^m$, with unit length $\| \bfb_i \|_2 = 1$, for $i=1,\ldots,k$.
	The low-rank factor $\bfZ = [\bfz_1, \ldots, \bfz_k]$ of T-LR-ADI is then given by the iteration:
	\begin{equation}\label{eq:Tang_ADI_real} 
		\begin{aligned}
		\bfz_i &= \sqrt{2\, s_i }  \left( \bfA - s_i\bfE \right)^{-1} \bfB_{\bot,i-1}\bfb_i, \\
		\bfB_{\bot,i} &= \bfB_{\bot,i-1} + \sqrt{2\, s_i }\bfE\bfz_i\bfb_i^T.
		\end{aligned}
	\end{equation}
\end{thm}
\begin{myproof}
	Due to Theorem~\ref{thm:ADI_Krylov}, it is sufficient to show that $\bfZ$---interpreted as a basis of a tangential Krylov subspace \eqref{eq:Krylov_t}---fulfills condition \eqref{eq:X}.
	As the re-formulated ADI iteration comprises a restart after every single step, it is actually sufficient to show this only for the first iterate $\bfz_1$.
	The matrices $\bfS$ and $\bfL$ in \eqref{eq:X} follow from the Sylvester equation \eqref{eq:syl}.
	We therefore rewrite the first equation of \eqref{eq:Tang_ADI_real} as $\left( \bfA - s_1\bfE \right) \bfz_1 \!=\! \sqrt{2\, s_1}  \bfB\bfb_1$, 
	from which we can identify $\bfS \!:=\! s_1$ and $\bfL \!:=\! \sqrt{2\, s_1}\bfb_1$.
	As by assumption $\bfb_1^H\bfb_1\!=\!1$, \eqref{eq:X} is solved by $\bfX \!:=\! 1$, and the approximate solution is $\Ph \!=\! \bfz_1\bfX^{-1}\bfz_1^T \!=\! \bfz_1\bfz_1^T$.
	It is left to prove that $\bfB_{\bot,1} \!=\! \bfB + \sqrt{2\, s_i }\bfE\bfz_1\bfb_1^T$.
	The residual factor is generally given by \eqref{eq:Bs_Krylov}.
	With $\bfV \!:=\! \bfz_1$, $\bfX \!:=\! 1$ and $\bfL \!:=\! \sqrt{2\, s_1}\bfb_1$, the result follows.
\end{myproof}

\subsection{Tangential ADI for Complex Conjugated Shifts}
Assume a complex shift $s_1 \in \mathbb{C}$ and tangential direction $\bfb_1 \in \mathbb{C}^m$, and compute the first iterate $\bfz_1$ by \eqref{eq:Tang_ADI_real} with $\sqrt{2\, \real (s_1) }$ instead of $\sqrt{2\, s_1}$.
In order to give a real approximation $\Ph \!=\! \bfZ\bfZ^H \in \mathbb{R}^{n \times n}$, the second iterate $\bfz_2$ has to be contained in $\colsp[\bfz_1, \overline{\bfz}_1]$.
However, the $\bfz_2$ resulting from the direct application of \eqref{eq:Tang_ADI_real} would not satisfy this.
For this reason, T-LR-ADI requires a slight modification for complex conjugated shifts.
\begin{thm}\label{thm:Tang_ADI_complex}
	Define $\bfB_{\bot,0} \!:=\! \bfB$, and assume for every iterate $i = 1,\ldots,k$ a complex shift $s_i \!\!\in\!\! \mathbb{C}$ with nonzero imaginary part, and a tangential direction $\bfb_i \!\!\in\!\! \mathbb{C}^m$, with $\| \bfb_i \|_2 \!\!=\!\! 1$.
	The low-rank factor $\bfZ \!\!=\!\! [\bfZ_1, \ldots, \bfZ_k]$ of T-LR-ADI, where each $\bfZ_i \in \mathbb{C}^{n \times 2}$ contains both ADI bases for the pair $(s_i,\bfb_i)$ and the complex conjugated pair $(\overline{s}_i,\overline{\bfb}_i)$, is then given by the iteration:
	\begin{equation}\label{eq:Tang_ADI_complex}
	\begin{aligned}
		\bfy_i &= \sqrt{2\, \real (s_i) } \left( \bfA - s_i\bfE \right)^{-1} \bfB_{\bot,i-1}\bfb_i, \\
		\alpha_i &= \bfb_i^H\overline{\bfb}_i \frac{\real(s_i)}{\overline{s}_i}, \quad
		\beta_i = \frac{1}{\sqrt{1- \alpha_i\overline{\alpha}_i}}, \\
		\bfZ_i &= \left[\bfy_i, \; \beta_i (\overline{\bfy}_i - \alpha_i\bfy_i)  \right], \\
		\bfB_{\bot,i} &= \bfB_{\bot,i-1} + \sqrt{2\, \real (s_i) }\bfE\bfZ_i\left[\bfb_i, \; \beta_i (\overline{\bfb}_i - \alpha_i\bfb_i)  \right]^H, 
	\end{aligned}
	\end{equation}
\end{thm}
\begin{myproof}
	First of all, note that $\beta_i$ always exists. 
	This is due to $\alpha_i\overline{\alpha}_i \!=\! \bfb_i^H\overline{\bfb}_i \bfb_i^T\bfb_i \cdot \real^2(s_i)\!/\!(s_i\overline{s}_i)$.
	Direct complex analysis shows that $0 \!\le\! \bfb_i^H\overline{\bfb}_i \bfb_i^T\bfb_i \!\le\! 1$, due to the assumption $\| \bfb_i \|_2 \!\!=\!\! 1$,
	and that $0 < \real^2(s_i)/(s_i\overline{s}_i) < 1$ because $s_i$ has nonzero imaginary part.
	Therefore: $0 < \alpha_i\overline{\alpha}_i < 1$ and $\beta_i$ always exists.
	
	The vector $\bfy_i$ points to the desired tangential direction and obviously, $\bfZ_i$ is a basis of $\colsp[\bfy_i, \overline{\bfy}_i]$.
	Similar to the proof of Theorem~\ref{thm:Tang_ADI_real}, it is sufficient to show, that $\bfZ_1$ fulfills condition \eqref{eq:X} and that $\bfB_{\bot,1}$ is the resulting low-rank factor of the residual.
	In order to construct $\bfZ_1$, we begin with $\bfV := [\bfy_1, \overline{\bfy}_1]$ as a basis of the desired tangential Krylov subspace.
	The Sylvester equation \eqref{eq:syl} for this $\bfV$ is satisfied by
	\begin{equation}
		\bfS_{\bfV} = \left[ \begin{array}{cc} s_1 & 0 \\ 0 & \overline{s}_1 \end{array} \right]
		\; \textrm{and} \;\;
		\bfL_{\bfV} = \sqrt{2\, \real (s_1) } \left[ \bfb_1, \, \overline{\bfb}_1 \right].
	\end{equation}
	Using $\bfb_1^H\bfb_1=1$, it can be readily verified that
	\begin{equation}
		\bfX_{\bfV} = \left[ \begin{array}{cc} 1 & \alpha_1 \\ \overline{\alpha}_1 & 1 \end{array} \right]
	\end{equation}
	solves the associated $2$-by-$2$ Lyapunov equation \eqref{eq:X}.
	The ADI approximation is given by $\Ph_1=\bfV\bfP_1\bfV^H$, with $\bfP_1 := \bfX_{\bfV}^{-1}$.
	Straightforward computation shows that $\bfP_1$ and its Cholesky factorization $\bfR\bfR^H := \bfP_1$ are given by
	\begin{equation}
		\bfP_1 = \bfX_{\bfV}^{-1} = \beta_1^2 \left[ \begin{array}{cc} 1 & -\alpha_1 \\ -\overline{\alpha}_1 & 1 \end{array} \right]
		\; \textrm{and} \;\;
		\bfR = \left[ \begin{array}{cc} 1 & -\alpha_1\beta_1 \\ 0 & \beta_1 \end{array} \right].
	\end{equation}
	Therefore, the low-rank Cholesky factorization of $\Ph_1\!=\!\bfV\bfR\bfR^H\bfV^H$ reads as $\Ph_1\!=\!\bfZ_1\bfZ_1^H$ with $\bfZ_1 \!=\! \bfV\bfR \!=\! \left[\bfy_1, \; \beta_1 (\overline{\bfy}_1 - \alpha_1\bfy_1)  \right]$.
	The associated residual factor $\bfB_{\bot,1}$ can be deduced from \eqref{eq:Bs_Krylov}:
	\begin{equation}
		\bfB_{\bot,1} = \bfB + \bfE\bfV\bfX_{\bfV}^{-1}\bfL_{\bfV}^H = \bfB + \bfE\bfZ_1(\bfL_{\bfV}\bfR)^H,
	\end{equation}
	with $\bfL = \bfL_{\bfV}\bfR = \sqrt{2\, \real (s_1) }\left[\bfb_1, \; \beta_1 (\overline{\bfb}_1 - \alpha_1\bfb_1)  \right]$, which completes the proof.
\end{myproof}
As every iterate $\bfZ_i$ is spanned by two complex conjugated vectors, the iteration \eqref{eq:Tang_ADI_complex} could also be formulated to deliver a real basis $\bfZ_i \in \mathbb{R}^{n \times 2}$, which was already done for the block ADI iteration \eqref{eq:BSADI}, see \cite{real_ADI} and \cite{Paed_real_Syl}.
The following theorem gives a possible real formulation of Theorem~\ref{thm:Tang_ADI_complex}.
\begin{thm}\label{thm:Tang_ADI_complex2}
	Define $\bfB_{\bot,0} \!:=\! \bfB$, and assume for every iterate $i = 1,\ldots,k$ a complex shift $s_i \!\!\in\!\! \mathbb{C}$ with nonzero imaginary part, and a tangential direction $\bfb_i \!\!\in\!\! \mathbb{C}^m$, with $\| \bfb_i \|_2 \!\!=\!\! 1$.
	The low-rank factor $\bfZ \!\!=\!\! [\bfZ_1, \ldots, \bfZ_k]$ of T-LR-ADI, where each real $\bfZ_i \in \mathbb{R}^{n \times 2}$ contains both ADI bases for the pair $(s_i,\bfb_i)$ and the complex conjugated pair $(\overline{s}_i,\overline{\bfb}_i)$, is then given by the iteration:
	\begin{equation}\label{eq:Tang_ADI_complex2}
	\begin{aligned}
		\bfy_i &= \sqrt{2\, \real (s_i) } \left( \bfA - s_i\bfE \right)^{-1} \bfB_{\bot,i-1}\bfb_i, \\
		\alpha_i &= \bfb_i^H\overline{\bfb}_i \frac{\real(s_i)}{\overline{s}_i}, \quad
		\beta_i = \frac{1}{\sqrt{1- \alpha_i\overline{\alpha}_i}}, \quad
		\gamma_i = \sqrt{1 + \real(\alpha_i)} \\
		\bfZ_i &= \frac{\sqrt{2}}{\gamma_i} \left[\real(\bfy_i), \; \beta_i \left(\imag(\alpha_i) \real(\bfy_i) + \gamma_i^2\imag(\bfy_i)\right)  \right], \\
		\bfB_{\bot,i} &= \bfB_{\bot,i-1} + \frac{2}{\gamma_i}\,\sqrt{\real (s_i) }\bfE\bfZ_i
		\left[\real(\bfb_i), \; \beta_i \left( \imag(\alpha_i)\real(\bfb_i) + \gamma_i^2\imag(\bfb_i)\right)  \right]^T, 
	\end{aligned}
	\end{equation}
\end{thm}
\begin{myproof}
	We start with the approximation of Theorem~\ref{thm:Tang_ADI_complex}, $\Ph_1=\bfV\bfP_1\bfV^H$, for which we show equivalence to $\Ph_1=\bfZ_1\bfZ_1^T$ of \eqref{eq:Tang_ADI_complex2}.
	The idea is to introduce an orthogonal basis transformation $\bfT \in \mathbb{C}^{2 \times 2}$:
	\begin{equation}\label{eq:T}
		\bfT := \frac{1}{\sqrt{2}}\left[\begin{array}{cc} 1 & -\imath \\ 1 & \imath \end{array}\right],
		\quad \bfT\bfT^H = \bfI.
	\end{equation}
	Then $\Ph_1\!=\!\bfV\bfP_1\bfV^H\!=\!\bfV\bfT\bfT^H\bfP_1\bfT\bfT^H\bfV^H\!=:\!\widetilde{\bfV}\widetilde{\bfP}_1\widetilde{\bfV}^H$, with $\widetilde{\bfP}_1\!:=\!\bfT^H\bfP_1\bfT \in \mathbb{R}^{2 \times 2}$ and $\widetilde{\bfV}\!:=\!\bfV\bfT \!=\! \sqrt{2}[\real(\bfy_1),\,\imag(\bfy_1)] \in \mathbb{R}^{n \times 2}$.
	Straightforward computation shows that $\widetilde{\bfP}_1$ and its Cholesky factorization $\widetilde{\bfR}\widetilde{\bfR}^H \!:=\! \widetilde{\bfP}_1$ are given by
	\begin{equation}
		\widetilde{\bfP}_1 = \bfT^H\bfP_1\bfT = \beta_1^2 \left[ \begin{array}{cc} 1-\real(\alpha_1) & \imag(\alpha_1) \\ \imag(\alpha_1) & 1+\real(\alpha_1) \end{array} \right]
		\; \textrm{and} \;\;
		\widetilde{\bfR} = \left[ \begin{array}{cc} \frac{1}{\gamma_1} & \frac{\beta_1\imag(\alpha_1)}{\gamma_1} \\ 0 & \beta_1\gamma_1 \end{array} \right].
	\end{equation}
	As $\Ph_1\!=\!\widetilde{\bfV}\widetilde{\bfR}\widetilde{\bfR}^H\widetilde{\bfV}^H$, the result for $\bfZ_1 \!=\! \widetilde{\bfV}\widetilde{\bfR}$ follows.
	The low-rank factor $\bfB_{\bot,1}$ of the residual is also derived from the proof of Theorem~\ref{thm:Tang_ADI_complex}:
	\begin{equation}
		\bfB_{\bot,1} = \bfB + \bfE\bfV\bfP_1\bfL_{\bfV}^H = 
		\bfB + \bfE\widetilde{\bfV}\widetilde{\bfP}_1\bfT^H\bfL_{\bfV}^H = 
		\bfB + \bfE\bfZ_1 (\bfL_{\bfV}\bfT\widetilde{\bfR})^H.
	\end{equation}
	With $\bfL_{\bfV}\bfT \!=\! 2\,\sqrt{\real(s_1)}[\real(\bfb_1),\,\imag(\bfb_1)] \in \mathbb{R}^{m \times 2}$, the result follows.
\end{myproof}
For a given set of shifts and tangential directions, both iterations of Theorem~\ref{thm:Tang_ADI_complex} and Theorem~\ref{thm:Tang_ADI_complex2} are equivalent, as they result in the same approximation $\Ph$ of the Lyapunov equation.
The benefit of the real formulation \eqref{eq:Tang_ADI_complex2} is that storage requirements and complex operations are reduced.

\subsection{Implementation}
Both ADI iterations for real \eqref{eq:Tang_ADI_real} and complex conjugated shifts, e.\,g. in the real version \eqref{eq:Tang_ADI_complex2}, can be combined.
An implementation of the resulting T-LR-ADI iteration is presented in Algorithm~\ref{alg:TangADI}.
The real low-rank factor $\bfZ=[\bfZ_1,\, \bfZ_2,\, ,\ldots,\, \bfZ_k]$ is computed iteratively, where $\bfZ_i$ has one column if $s_i$ is real, and two columns if $s_i$ is complex.
Additionally, the low-rank factor $\Bs$ of the residual is iteratively computed.
\renewcommand{\algorithmicrequire}{\textbf{Input:}}
\renewcommand{\algorithmicensure}{\textbf{Output:}}
\begin{algorithm}[!ht]\caption{Tangential-Low-Rank-ADI (T-LR-ADI)}\label{alg:TangADI}
\begin{algorithmic}[1]
		\Require $\bfE$, $\bfA$, $\bfB$, $tol$
		\Ensure Approximation $\Ph = \bfZ\bfZ^T$ and residual $\bfR = \Bs\Bs^T$
		\State initial choice of $s_1 \in \mathbb{C}$ and $\bfb_1 \in \mathbb{C}^m$ with $\| \bfb_1 \|_2 = 1$
		\State {$\bfZ=[\;]$, $\Bs = \bfB$}
		\Repeat  
		\State solve $\left( \bfA - s_i\bfE \right)\bfy=\Bs\bfb_i$ for $\bfy$
		\If{$s_i \in \mathbb{R}$ }
			\State $\bfZ_i = \sqrt{2s_i}\bfy$, $\; \bfL_i = \sqrt{2s_i}\bfb_i$
		\Else 
			\State $\alpha = \bfb_i^H\overline{\bfb}_i \frac{\real(s_i)}{\overline{s}_i}, \quad \beta = \frac{1}{\sqrt{1- \alpha\overline{\alpha}}}, \quad \gamma = \sqrt{1 + \real(\alpha)}$
			\State $\bfZ_i = \frac{2}{\gamma}\,\sqrt{\real(s_i)} \left[\real(\bfy), \; \beta \left( \imag(\alpha)\real(\bfy) + \gamma^2\imag(\bfy)\right)  \right]$
			\State $\bfL_i = \frac{2}{\gamma}\,\sqrt{\real(s_i)}\left[\real(\bfb_i), \; \beta \left( \imag(\alpha) \real(\bfb_i) + \gamma^2\imag(\bfb_i)\right)  \right]^T$
		\EndIf
		\State $\bfZ = [\bfZ,\bfZ_i]$
		\State $\Bs = \Bs + \bfE\bfZ_i\bfL_i^T$
		\State determine $s_{i+1}$ and $\bfb_{i+1}$ with $\| \bfb_{i+1} \|_2 = 1$
		\Until{$\| \bfR \|_2 \!=\! \max\,\Lambda\!\left(\Bs^T\Bs\right)  < tol \, \| \bfB^T\bfB \|_2$}
\end{algorithmic}
\end{algorithm}

\subsection{Application Aspects}
T-LR-ADI by Algorithm~\ref{alg:TangADI} represents a generalization of the block ADI iteration \eqref{eq:BSADI}.
To demonstrate this, assume $m$ equal real shifts $s_1 \!\!=\!\! s_2 \!\!=\!\! \ldots \!\!=\!\! s_m$ and an orthonormal basis $\{\bfb_1,\ldots,\bfb_m\}$ as tangential directions.
Then the basis $[\bfz_1,\ldots,\bfz_m]$ of T-LR-ADI is equal to the first iterate of the block ADI iteration \eqref{eq:LRADI} or \eqref{eq:BSADI}.
Therefore, the block ADI iteration is a special case of T-LR-ADI, where the latter provides an additional degree of freedom:
instead of being restricted to the whole block, we can pick only certain directions of our choice.
As the numerical experiments will show, neglecting non-relevant directions in the block can even improve approximation quality.

The typical employment of the block ADI iteration \eqref{eq:BSADI} is recycling an a priori chosen set of shifts $\mathcal{S} = \{s_1,\ldots,s_k\}$.
The set $\mathcal{S}$ is usually found by an optimization on the spectrum of $\bfE^{-1}\bfA$, see \cite{diss_saak} for an overview.
However, such a choice has two drawbacks for application with T-LR-ADI.

On the one hand, the optimization is solely performed on the spectrum of $\bfE^{-1}\!\bfA$ and it neglects the effect of the right-hand side $\bfB\bfB^T$ on the Lyapunov equation.
However, it is crucial to also take $\bfB$ into account to reasonably select tangential directions $\bfb_i$.

On the other hand, in order to benefit from the new degree of freedom in T-RL-ADI, each shift should only be reused $r\!<\!m$ times (together with the orthonormal basis $\{\bfb_1,\ldots,\bfb_r\}$)---because otherwise the block version could have been directly applied.
It would be hard to a priori guess a suitable set of shifts $\mathcal{S}$, that certainly leads to convergence after recycling each shift only $r<m$ times.
Therefore, instead of selecting the set $\mathcal{S}$ before the ADI iteration, it is more reasonable to successively determine a new shift and ($r$) tangential direction(s) for the next iterate.
We believe that such a selection strategy, that is based on the currently available approximation (this could be the recent ADI basis $\bfZ_i$ or the current residual factor $\bfB_{\bot,i}$) can outperform a choice of shifts disconnected from the growing ADI approximation.

However, this represents a research direction in its own right.
We only propose a first selection strategy in Section~\ref{sec:adapt}, by using basic ideas that are already available in the literature.
The scope of this work is rather to disclose the potential of T-LR-ADI, which will be revealed by the numerical examples in Section~\ref{sec:num_ex}.

\section{Adaptive Choice of Shifts and Tangential Directions}\label{sec:adapt}
The eigenvalues of the projected matrices \eqref{eq:ArErBr}, associated with the ADI iteration, are the mirrored shifts $s_i$.
This reveals the strong connection of the ADI iteration to $\H2$ optimal model order reduction, from which we borrow the shift selection strategy.
Gugercin et al.\,\cite{H2_gugercin} proposed an \emph{iterative rational Krylov algorithm} (IRKA) that, upon convergence, yields a locally $\H2$ optimal reduced model as follows:
for a given set of shifts, the matrices $\bfA$ and $\bfE$ are projected onto the (tangential) rational Krylov subspace like in \eqref{eq:ArErBr};
subsequently, a new (tangential) rational Krylov subspace is computed, using as shifts the mirrored Ritz values with respect to the imaginary axis
(and as tangential directions the rows of that $\bfB_k$ that corresponds to diagonalized $\bfA_k$, $\bfE_k$);
this procedure is iterated until convergence occurs.

The basic idea is to apply IRKA in every iteration of T-LR-ADI.
The main difference is that it is not necessarily performed until convergence occurs:
instead, we predefine a maximum number of iterations $n_{max}$, that can even be as small as $1$.
As we are interested in only one real or two complex conjugated shifts, the reduced order of IRKA in each iteration is either $1$ or $2$.
This should usually lead to fast convergence of IRKA, because typically the convergence speed of IRKA deteriorates the larger the reduced order is.
As initial projection matrix in IRKA we use the current ADI basis $\bfZ_i \!\!=:\!\! \bfV$ at step $14$ in Algorithm~\ref{alg:TangADI} together with $\bfW \!\!:=\!\! \bfV$.
From the projected matrices \eqref{eq:ArErBr}, new shifts and tangential directions are determined as described above.
After at most $n_{max}$ iterations of IRKA, a new (complex) shift and associated tangential direction are found for the next iterate in T-LR-ADI.

A possible implementation is sketched in Algorithm~\ref{alg:IRKA_update}, where $\bfM(1,:)$ denotes the first row of a matrix $\bfM$.
\begin{algorithm}[!ht]\caption{update scheme for T-LR-ADI}\label{alg:IRKA_update}
\begin{algorithmic}[1]
		\Require $\bfE$, $\bfA$, $\Bs$, $\bfZ_i$, $tol$, $n_{max}$
		\Ensure $s_{i+1}$ and $\bfb_{i+1}$ with $\| \bfb_{i+1} \|_2 = 1$
		\State $\bfV = \bfZ_i$
		\For{$j = 1,\ldots,n_{max}$}  
			\State $\bfA_j \!=\! (\bfV^H\bfE\bfV)^{-1}\bfV^H\bfA\bfV$, $\; \bfB_j \!=\!(\bfV^H\bfE\bfV)^{-1} \bfV^H\Bs$
			\State compute $\Lambda(\bfA_j)$: $\bfA_j = \bfU_j\bfD_j\bfU_j^{-1}$, $\widetilde{\bfB}_j = \bfU_j^{-1}\bfB_j$
			\If{relative change of $\Lambda(\bfA_j) < tol$ } break;
			\EndIf
			\State solve $\bfA\bfV + \bfE\bfV\bfD_j = \Bs\widetilde{\bfB}_j^H$ for $\bfV$ 
			{\hfill // compute basis of tangential Krylov subspace}
		\EndFor
		\State $s_{i+1} = -\bfD_j(1,1)$, $\bfb_{i+1} = \widetilde{\bfB}_j(1,:)/\|\widetilde{\bfB}_j(1,:) \|_2$
		\If {$s_{i+1} \in \mathbb{R}$ and $\bfV \in \mathbb{R}^{n \times 2}$}
			\State $s_{i+2} = -\bfD_j(2,2)$, $\bfb_{i+2} = \widetilde{\bfB}_j(2,:)/\|\widetilde{\bfB}_j(2,:) \|_2$
		\EndIf
\end{algorithmic}
\end{algorithm}
The Sylvester equation at step $7$ should be merely interpreted as a notation that allows to include both reduced orders $1$ and $2$:
its solution $\bfV$ is a basis of the tangential rational Krylov subspace with the mirrored Ritz values as shifts and the rows of $\widetilde{\bfB}_j$ as tangential directions;
therefore, one (complex) linear system solve is required in step~$7$;
only if both eigenvalues of $\bfA_j$ happen to be real, two system solves are required.
Note that using Algorithm~\ref{alg:IRKA_update} at step $14$ of T-LR-ADI, the solve at step $4$ of the latter becomes redundant, as $\bfy$ is already contained in the final $\bfV$ of Algorithm~\ref{alg:IRKA_update}.
After $n_{max}$ iterations of Algorithm~\ref{alg:IRKA_update}, two special cases can occur:
if $\bfD_j(1,1)$ or $\bfD_j(2,2)$ is contained in the right-half of the complex plane, we simply take $s_{i+1} = \bfD_j(1,1)$ and $s_{i+2} = \bfD_j(2,2)$, respectively;
if $\bfA_j$ has two real eigenvalues, we apply both $s_{i+1}$ and $s_{i+2}$ to T-LR-ADI and initialize Algorithm~\ref{alg:IRKA_update} the next time with $\bfV:=[\bfZ_{i+1},\bfZ_{i+2}]\in\mathbb{R}^{n \times 2}$.

The main numerical effort in Algorithm~\ref{alg:IRKA_update} is the computation of the tangential Krylov subspace in step $7$, which amounts to the solution of one (complex) $n$-by-$n$ linear system per iteration (only if $\bfA_j$ has two real eigenvalues, two real solves are required).
If we choose $n_{max} = 1$, every solution of a linear $n$-by-$n$ system is added to the low-rank factor $\bfZ$, and Algorithm~\ref{alg:IRKA_update} does not increase the overall numerical effort of T-LR-ADI---besides a $2$-by-$2$ eigenvalue problem.
Only if $n_{max} > 1$, solutions are discarded.
Therefore, the maximum number of iterations $n_{max}$ should be chosen as small as possible for low numerical effort.
But $n_{max}$ should also be chosen as large as necessary, because discarding linear system solves typically increases the quality of shifts.
This causes a trade-off, because the quality of shifts affects the convergence speed of T-LR-ADI and the dimension of the final ADI basis $\bfZ$.
A general rule for selecting $n_{max}$ is challenging, and probably impossible, as first numerical tests suggest that this is strongly dependent on the problem:
on the one hand, for problems with excellent convergence behavior, $n_{max} = 1$ would not only lead to the smallest overall execution time of T-LR-ADI, but also the number of columns in the final $\bfZ$ would be only marginally larger than for high $n_{max}$ yielding better shifts.
On the other hand, for problems with poor convergence behavior, T-LR-ADI might not converge at all in reasonable time for small $n_{max}$ in Algorithm~\ref{alg:IRKA_update}, making a higher choice of $n_{max}$ essential.
We compared different values of $n_{max}$ in our numerical examples.

A related adaptive shift selection strategy was suggested in \cite{Paed_ADI_Shifts}:
for a given set of initial shifts $s_1,\ldots,s_r$, $r\le m$, the block ADI iteration \eqref{eq:BSADI} is performed.
The last iterate $\bfZ_r \!\in\! \mathbb{C}^{n \times m}$ is used to define $\bfV\!\!:=\!\! [\real(\bfZ_r),\textrm{Im}(\bfZ_r)]$.
Then a projection is performed like in step~$3$ of Algorithm~\ref{alg:IRKA_update}.
From the subsequent spectral decomposition, $2m$ Ritz values follow.
Then the $r\le 2m$ Ritz values are selected, which are contained in the left half of the complex plane.
Their mirror images are then used as new shifts in the block ADI iteration.
However, it is unclear, in which order the shifts $s_1,\ldots,s_{r}$ have to be used, as this generally affects the Ritz values after projection.

Both adaptive shift selection strategies are designed for problems with $m\!>\!1$, however, they can be readily applied to problems with $m\!=\!1$.
For T-RL-ADI then simply set all $\bfb_i\!=\!1$.
If additionally $n_{max}$ in Algorithm~\ref{alg:IRKA_update} was set to $n_{max}\!=\!1$, both strategies would result in the same sequence of shifts.

The initial shifts suggested in \cite{Paed_ADI_Shifts} are the mirrored Ritz values after projection with $\bfB$.
Due to the mirroring of shifts and eigenvalues, we propose to choose one (or a few) mirrored eigenvalues of $\Lambda(\bfE^{-1}\bfA)$ instead.
For example, the ($m$) eigenvalue(s) with smallest magnitude can usually be computed with low numerical effort.
A possible implementation to compute $s_1$ and also $\bfb_1$ in MATLAB is given in the following:
\begin{align*}
	& \texttt{[v,s1] = eigs(A,E,1,'SM'); s1 = -s1;} \\
	& \texttt{b1 = v'*(E} \text{\textbackslash} \texttt{B); b1 = b1'/norm(b1);}
\end{align*}

\section{Numerical Example}\label{sec:num_ex}
We do not perform extensive numerical tests here, as this is out of the scope of the paper.
The purpose of the numerical examples is rather to emphasize that T-LR-ADI can be a suitable alternative to the block ADI iteration---depending on the problem.

The first example we consider is a semi-discretized heat transfer problem for optimal cooling of steel profiles from the Oberwolfach model reduction benchmark collection\footnote{Available at http://portal.uni-freiburg.de/imteksimulation/downloads/benchmark.}
with $n\!\!=\!\!1,\!357$ and $m\!\!=\!\!7$.
We performed T-LR-ADI with the initial choice of $s_1$ and $\bfb_1$ as suggested above.
In this first example, we would like to disclose the efficiency of the tangential ADI basis in approximating $\bfP$---when compared to the block basis with equal set of shifts.
For the time being, we are therefore not concerned with execution time, and choose a (practically inappropriate) $n_{max}\!=\!200$, to find preferably optimal shifts.
As the problem has a real spectrum, the reduced order for the IRKA approach in Algorithm~\ref{alg:IRKA_update} was one and only real shifts and tangential directions were adaptively computed.
The resulting sequence of shifts was subsequently applied to the block version of the ADI iteration \eqref{eq:BSADI}.
Figure~\ref{fig:Vergleich_tan_block} shows the relative error of the approximation $\Ph \!=\! \bfZ\bfZ^H$ after each iteration, where the full-rank $\bfP$ was calculated with the command \texttt{lyapchol} in MATLAB.
For $i$ iterations, the low-rank factor $\bfZ$ has $i$ columns for T-LR-ADI and $mi\!=\!7i$ columns for the block ADI iteration.
After $194$ iterations for instance, the low-rank factor $\bfZ$ of the block ADI iteration reaches $194\!\cdot\! 7 \!=\! 1,\!358 \!>\! n$ columns with the relative error $5.3\!\cdot\! 10^{-5}$; whereas $\bfZ$ of T-LR-ADI has only $194$ columns with the relative error $9.7\cdot 10^{-10}$.
Therefore, compared to the block ADI iteration with the same set of shifts, T-LR-ADI can yield a better approximation with several orders of magnitude---despite a significantly smaller rank.
\newlength{\maxwidth}
\setlength{\maxwidth}{\columnwidth}
\newlength{\myheight}
\setlength{\myheight}{0.382\columnwidth}
\begin{figure}[!htb]
	\footnotesize
	\centering
	\includegraphics[width=\maxwidth,height=\myheight]{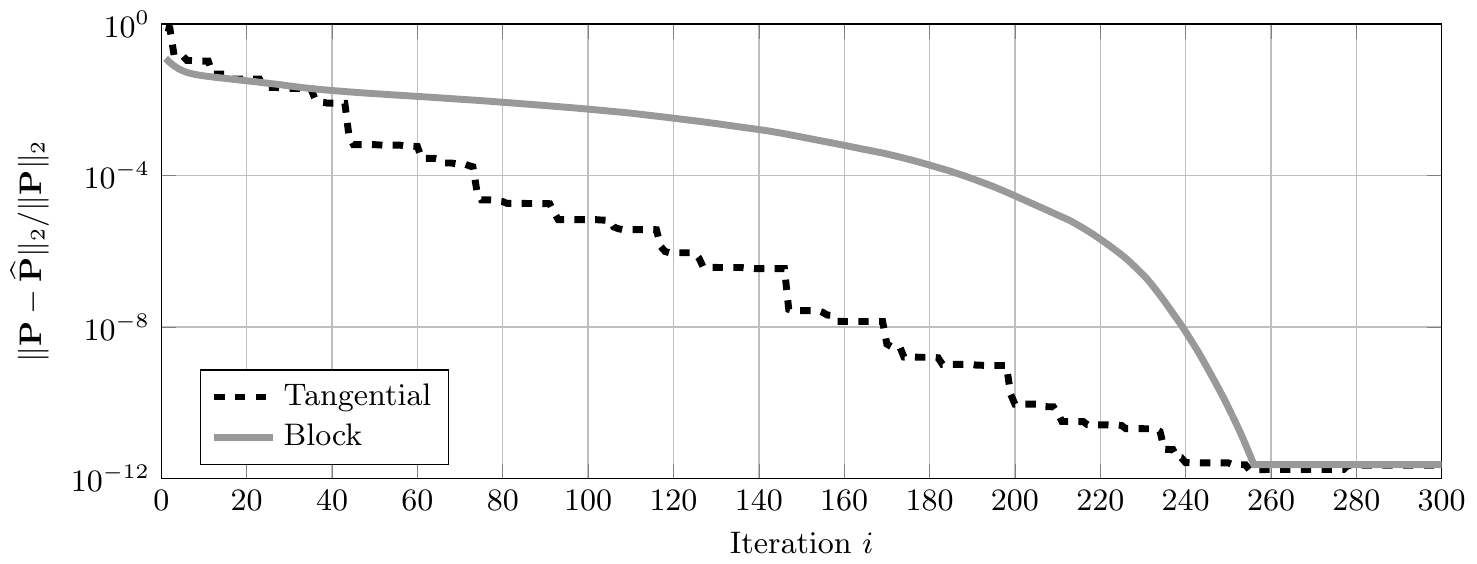}	
	\caption{Relative error for steel profile, $n = 1,\!357$, $m = 7$}
	\label{fig:Vergleich_tan_block}
\end{figure}

Typically, the goal of ADI is to compute an approximation as fast as possible.
In a block iteration, the $n$-by-$n$ linear system is solved with a multiple right-hand side, whereas for every column in T-LR-ADI, one linear system with a new coefficient matrix has to be solved.
Therefore, the numerical effort per column in $\bfZ$ is smaller for a block iteration.
The question is now, if T-LR-ADI can overcome this drawback, by requiring less columns for convergence (as Figure~\ref{fig:Vergleich_tan_block} suggests).
To answer this question in the following, we compare T-RL-ADI with different choices of $n_{max}$ and denoted as T-LR-ADI($n_{max}$), to the adaptive block iteration, denoted as B-LR-ADI and given in \cite{Paed_ADI_Shifts} (it also computes a real ADI basis for complex conjugated shifts).
There, B-LR-ADI(a) refers to the originally suggested choice of initial shifts (mirrored Ritz values after projection with $\bfB$) and B-LR-ADI(b) denotes the here proposed initial selection ($m$ mirrored eigenvalues with smallest magnitude).

We used the steel model with $n\!=\!5,\!177$ as an example, and ran all methods until the relative residual error $\|\bfR\|_2$/$\|\bfB\bfB^T\|_2$ dropped below $10^{-12}$.
Table~\ref{tab:perf_steel} shows the results, where $t_{\textrm{shift}}$ denotes the time for computing the initial shifts (and tangential direction), $t_{\textrm{ADI}}$ the time for the ADI iteration and $t_{\textrm{total}}$ the respective sum.
All times are stated in seconds.
The total number of iterations is denoted by $k$ and $dim$ states the corresponding number of columns in $\bfZ$.
\begin{table}[htb]
	\centering
	\caption{Steel profile, $n \!=\! 5,\!177$, $m \!=\! 7$}
		\begin{tabular}{lccccc}
			& $t_{\textrm{shift}}$ & $t_{\textrm{ADI}}$ & $t_{\textrm{total}}$ & $k$ & $dim$ \\
			\hline
				     T-LR-ADI($1$) & 2.9e-1 & 8.71 & 9.00 & 367 & 367 \\
				     T-LR-ADI($2$) & 2.9e-1 & 14.5 & 14.8 & 346 & 346 \\
				     T-LR-ADI($4$) & 2.9e-1 & 26.6 & 26.9 & 333 & 333 \\
				     T-LR-ADI($8$) & 2.9e-1 & 49.4 & 49.7 & 322 & 322 \\
				    T-LR-ADI($50$) & 2.9e-1 & 192  & 192  & 313 & 313 \\
					  B-LR-ADI(a)  & 1.3e-3 & 5.14 & 5.14 &  71 & 497 \\
					  B-LR-ADI(b)  & 1.9e-1 & 3.91 & 4.10 &  58 & 406 \\
			\hline
		\end{tabular}
	\label{tab:perf_steel}
\end{table}

It can be seen, that the new initial choice of shifts B-LR-ADI(b) outperforms the originally suggested B-LR-ADI(a) both in terms of size of $\bfZ$ and execution time.
As expected, for all values of $n_{max}$, T-LR-ADI leads to a (much) smaller low-rank factor $\bfZ$, compared to B-LR-ADI.
However, as the problem exhibits good convergence behavior, B-LR-ADI converges faster and the smaller basis $\bfZ$ of T-LR-ADI does not pay off here.

\begin{figure}[!htb]
	\footnotesize
	\centering
	\includegraphics[width=\maxwidth,height=\myheight]{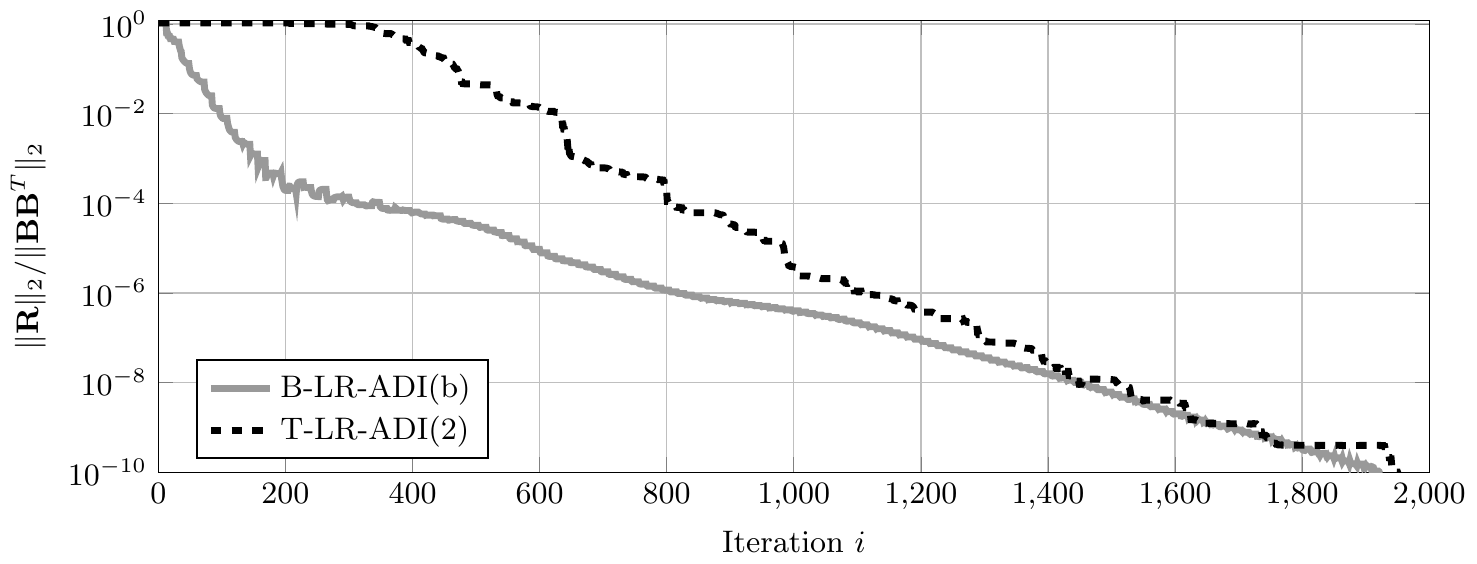}	
	\caption{Relative residual norm for butterfly gyro, $n = 34,\!722$, $m = 12$}
	\label{fig:Vergleich_Gyro}
\end{figure}

We expect problems with complex spectra to have worse convergence behavior.
As an example we consider the \emph{butterfly gyro} model with $n\!=\!34,\!722$ from the same benchmark collection and solve the dual Lyapunov equation with $\bfA^T$, $\bfE^T$ and $\bfC^T$ instead of $\bfA$, $\bfE$ and $\bfB$, respectively; then $m\!=\!12$.
Calculations were carried out in MATLAB on a dual-core processor with 2.27 GHz and 48 GB RAM; the linear systems were solved with the backslash operator.
The initial choice B-LR-ADI(a) resulted in $m$ zeros, so we could only employ B-LR-ADI(b) here.
We ran it together with T-LR-ADI(2) until the relative residual error $\|\bfR\|_2$/$\|\bfB\bfB^T\|_2$ dropped below $10^{-10}$.
T-LR-ADI(2) took $3.4$ hours for $1,\!951$ iterations, resulting in $\bfZ$ with $3,\!288$ columns (i.\,e. both real and complex conjugated shifts were found), 
and B-LR-ADI(b) took $4.8$ hours for $1,\!922$ iterations, yielding an ADI basis with $23,\!064$ columns.
Figure~\ref{fig:Vergleich_Gyro} shows that for B-LR-ADI(b) the relative residual error rapidly drops in the beginning.
This is due to the high $m\!=\!12$, which leads to a fast accumulation of columns in the block iteration.
However, its basis $\bfZ$ also incorporates dispensable directions that slow down convergence after approximately $200$ iterations.
In contrast, the relative residual error does not at all decrease in the beginning of T-LR-ADI.
However, after a sufficient amount of columns in $\bfZ$, convergence speeds up.
And in the end, T-LR-ADI outperforms B-LR-ADI in this example both in terms of final rank and execution time.

Figure~\ref{fig:Vergleich_Gyro} suggests that it should be possible to combine the advantages of both approaches:
it sounds reasonable to start with the block ADI iteration---to rapidly accumulate a preparative basis---and then switch to the tangential version---to limit the rank of the final approximation.
A more sophisticated approach would assign an optimal number of $r$, with $1\le r\le m$, tangential directions for every shift $s_i$.
Then computation time could be saved---similar to the block iteration---by solving linear systems with multiple right-hand sides in T-LR-ADI.
This would allow to adaptively ``blend'' between the tangential approach and the block method.
The potential of such a \emph{multi directional} tangential method was recently unveiled in the context of RKSM \cite{Druskin_adapt_MIMO}---which indicates high potential also for the ADI iteration.
However, this requires enhanced selection strategies, which is a research topic in its own right.
The purpose of the present work is only to introduce tangential directions in the ADI iteration and to present a first approach to adaptively select shifts and tangential directions.

\section{Conclusion}
A new ADI iteration for the solution of large-scale Lyapunov equations was introduced, which allows to add a new degree of freedom: the application of tangential directions.
Together with an adaptive shift selection strategy, the new tangential method leads to approximations with less dimension and storage requirements than existing approaches based on blocks.
Depending on the problem, tangential ADI can outperform existing strategies additionally with respect to execution time.

Future work includes optimized strategies for selecting shifts and tangential directions and the generalization to large-scale Sylvester equations.

\bibliographystyle{plainnat}
\bibliography{ref}

\end{document}